\documentclass[12pt,draft]{amsart}
\usepackage{amssymb, amsmath}
\usepackage{url}
\newtheorem{theorem}{Theorem}[section]

\theoremstyle{remark}

\theoremstyle{definition}

\newtheorem{example}[theorem]{Example}

\numberwithin{equation}{section}
\numberwithin{theorem}{section}

\let\abs=\envert

\newcommand{\N}{{\mathbb N}}
\newcommand{\R}{{\mathbb R}}
\newcommand{\C}{{\mathbb C}}

\newcommand{\fn}{\!:\!}

\providecommand{\abs}[1]{\lvert#1\rvert}


\begin{document}

\date{Preprint January 14, 2026.}

\title{Summing series using recurrence relations}
\author{Erik Talvila}
\address{Department of Mathematics \& Statistics\\
University of the Fraser Valley\\
Abbotsford, BC Canada V2S 7M8}
\email{Erik.Talvila@ufv.ca}

\begin{abstract}
Power series in which the summand satisfies a linear recurrence relation
with polynomial coefficients are shown to be the solution of a linear
differential or algebraic equation.  Solving the associated 
differential or algebraic equation yields a closed form for the series.  This method is used to
sum several series 
and to solve two {\it American Mathematical Monthly} problems.
\end{abstract}

\maketitle

\section{Introduction}\label{sectionintroduction}
If we have a series we wish to sum, that is, write in closed form, it may happen
that the summand has useful properties.  If the summand satisfies a 
linear recurrence
relation with polynomial coefficients then the associated power 
series (known as the generating function) will be the solution of a linear
algebraic or differential equation. 

We begin with
a simple example to illustrate this method, which involves solving an algebraic
equation (Section~\ref{sectioncosine}).
Then, in Section~\ref{sectiontheory}, a small amount of theory will be
presented and some indications of the literature on the subject.  Next, another
algebraic example (Section~\ref{section10977}).  The paper
then concludes with three examples that involve solving a differential
equation (Sections~\ref{sectionbinomial}, \ref{sectionbessel}, \ref{section11274}).
Along the way we will encounter the binomial theorem, a generating function
for Bessel functions and two {\it American Mathematical Monthly} problems.

\section{Cosine example}\label{sectioncosine}
\begin{example}
Suppose we wish to obtain a closed form expression for $\sum_{n=0}^N\cos(n)$.
The cosine addition formula leads to the identity
$\cos(n+1)-2\cos(1)\cos(n)+\cos(n-1)=0$.
Now bring in the generating function.
Let $S_N(x)=\sum_{n=0}^Nx^n\cos(n)$. Multiply the cosine recurrence relation by
$x^n$ and sum from $n=0$ to $n=N$.  Re-index the first and third sum,
\begin{align*}
0&=\sum_{n=0}^Nx^n\cos(n+1)-2\cos(1)\sum_{n=0}^Nx^n\cos(n)+\sum_{n=0}^Nx^n\cos(n-1)\\
&=\sum_{n=1}^{N+1}x^{n-1}\cos(n)
-2\cos(1)\sum_{n=0}^Nx^n\cos(n)+\sum_{n=-1}^{N-1}x^{n+1}\cos(n)\\
&=\frac{S_N(x)}{x}+x^N\cos(N+1)-\frac{1}{x}-2\cos(1)S_N(x)+xS_N(x)-x^{N+1}\cos(N)
+\cos(1).
\end{align*}
Solving this algebraic equation for $S_N(x)$ gives
$$
S_N(x)=\frac{x^{N+2}\cos(N)-x^{N+1}\cos(N+1)-x\cos(1)+1}{x^2-2\cos(1)x+1}.
$$
Putting $x=1$ gives 
$$
\sum_{n=0}^N\cos(n)
=\frac{\cos(N)-\cos(N+1)-\cos(1)+1}{2(1-\cos(1))}.
$$
If $\abs{x}<1$ then
$$
\sum_{n=0}^\infty x^n\cos(n)=\frac{1-x\cos(1)}{1-2\cos(1)x+x^2}.
$$

This sum can also be obtained by taking the real part of the
geometric series $$
\sum_{n=0}^N\left({\rm sgn}(x)e^{\log\abs{x}+i}\right)^n.
$$
There are many examples of this type in \cite{wilf}.
We will see examples below that cannot be summed in this manner.
\end{example}

\section{Recurrence relations with polynomial coefficients}\label{sectiontheory}
Our method of summation is based on the identities
\begin{equation}
\sum_{n=0}^Na_{n+1}x^n=\sum_{n=1}^{N+1}a_{n}x^{n-1} =x^{-1}\sum_{n=0}^Na_{n}x^n
+a_{N+1}x^{N}-a_0x^{-1},\label{constant}
\end{equation}
which was used in the cosine example, and
\begin{equation}
\sum_{n=0}^Nna_{n}x^n=x\sum_{n=1}^Nna_{n}x^{n-1}=x\frac{d}{dx}\sum_{n=0}^Na_{n}x^n,\label{polynomial}
\end{equation}
which will be used in Sections~\ref{sectionbinomial}, \ref{sectionbessel} and 
\ref{section11274} below.

In the cosine example, the summand satisfied a linear recurrence relation
with constant coefficients.  Let's see what happens when the coefficients
of the recurrence relation are polynomials.  Suppose we have a sum
$\sum_{n=1}^Na_n$ and the recurrence relation
$$
\alpha_{n,k}a_{n+k}+\alpha_{n,k-1}a_{n+k-1}+\ldots+\alpha_{n,0}a_{n}=0
$$
where the $\alpha$ coefficients are polynomials in $n$.  Consider the
term $\alpha_{n,m}a_{n+m}$ for some $0\leq m\leq k$.  This will be 
multiplied by $x^n$ and then re-indexed to become $\alpha_{n-m,m}a_{n}x^{n-m}$.
If $\alpha_{\,\,\cdot,m}$ is a linear function, say, $\alpha_{n-m,m}=p(n-m)+q$
for constants $p$ and $q$,
then we can write $\alpha_{n-m,m}a_{n}x^{n-m}=p(n-m)a_nx^{n-m}+qa_nx^{n-m}
=pxd/dx[a_{n}x^{n-m}]+qa_nx^{n-m}$.  This term will then contribute a
term $xS_N'(x)$ and a term $S_N(x)$.  Similarly with the other terms in the
recurrence relation.  For the series $S_N$ we would then get a first order
linear differential equation in $S_N$ with linear coefficients.  

By the
same re-indexing and identification of derivatives, if the coefficients of
the recurrence relation were polynomials of degree $r$ then the sum $S_N$
would satisfy an $rth$ order differential equation.  Differential equations
with polynomial coefficients are typically solved using series methods.
That's exactly what we don't want here.  So, the case of interest is when
the recurrence relation has linear coefficients.  Then the first order
differential equation can always be solved in terms of an integral and
this will give a closed form expression for our series, perhaps in terms
of an unevaluated integral.

Of course there are many ways of summing series, such as evaluating Taylor
or Fourier series at a point, use of residues for analytic functions,
using Laplace transforms, employing
transformations and rearrangements, plus numerous specialized methods
\cite{henrici},
\cite{kalman},
\cite{knopp}, \cite{wheelon1954}, \cite{wheelon1955}.

Given $F(n,k)$,
Gosper's algorithm produces $G$ such that $F(n+1,k)-F(n,k)=G(n,k+1)-G(n,k)$
or says no such $G$ exists.
Upon summation on $k$ the right side telescopes and this shows $\sum_{k=0}^N F(n,k)$
is independent of $n$.  Substituting a particular value of $n$ gives this sum.
The algorithm is applicable when $F(n,k+1)/F(n,k)$ is a rational function of $k$.
See \cite{brouwer}, \cite{nemes} and
\cite{petkovsekwilfzeilberger}, where various generalisations
are also described.  
These methods have
been implemented in computer algebra systems such as Maple and
Mathematica.  Extension to differential
equations for power series is given in \cite{koepf}.

We'll solve an {\it American Mathematical Monthly} problem with a 
recurrence relation that has an algebraic
solution.  Then the remaining examples will involve differential equations.

\section{
American Mathematical Monthly Problem 10977}\label{section10977}
{\it
Let $\beta$ and $t$ be real numbers with $t>0$. Let ${\rm erf}(x)=
\frac{2}{\sqrt{\pi}}
\int_0^xe^{-t^2}\,dt$, let $\gamma(x,y)=\int_0^y t^{x-1}e^{-t}\,dt$, and 
let $\Gamma(x)=\int_0^\infty t^{x-1}e^{-t}\,dt$.  Show that
$$
\sum_{n=0}^\infty\left(\frac{\beta^2}{1+\beta^2}\right)^n\frac{\gamma(n+1/2,
(1+\beta^2)t)}{\Gamma(n+1/2)}=(1+\beta^2){\rm erf}(\sqrt{(1+\beta^2)t})
-e^{-t}\beta\sqrt{1+\beta^2}\,{\rm erf}(\beta\sqrt{t}).
$$
}
See \cite{10977}.\\

\noindent {\it Solution.}
It is assumed $x>0$.

Integrate by parts to get
\begin{equation}
\gamma(x+1,y)    =  \int_{t=0}^y t^x e^{-t}\,dt
  =  -y^x e^{-y} +x\int_{t=0}^y t^{x-1} e^{-t}\,dt
  = x\,\gamma(x,y)-y^x e^{-y}.\label{1}
\end{equation}
The following identities follow from the given definitions:
\begin{eqnarray}
\Gamma(x+1) & = & x\,\Gamma(x)\label{2}\\
\Gamma(1/2) & = & \sqrt{\pi}\label{3}\\
\gamma(1/2, x^2) & = & \sqrt{\pi}\,\text{\rm erf}(x)\label{4}\\
\Gamma(n+1/2) & = & \sqrt{\pi}\,\frac{1\cdot 3\cdot5\cdots (2n-1)}{2^n}.\label{5}
\end{eqnarray}
We will also use the identity
\begin{equation}
\int_0^{\pi/2}\cos^{2m+1}\!\theta\,d\theta=\frac{2^m m!}
{1\cdot 3\cdot5\cdots (2m+1)}.\label{6}
\end{equation}

Let $a_n(y)=\gamma(n+1/2,y)/\Gamma(n+1/2)$ and $b_n(y)=y^{n+1/2}e^{-y}/
\Gamma(n+3/2)$.  Using \eqref{1} and \eqref{2}, we have
$a_{n+1}(y)=a_n(y)-b_n(y)$.

Define $T_y(x)=\sum_{n=0}^\infty a_n(y)x^n$.  Since $\abs{\gamma(x,y)}\leq
\Gamma(x)$, this series converges absolutely for $\abs{x}<1$ and all $y\geq 0$.
Due to the recurrence relation for $a_n$ the series can be summed in
closed form because it telescopes:
\begin{eqnarray*}
\sum_{n=0}^\infty a_{n+1}(y) \,x^{n+1} & = &
\sum_{n=0}^\infty a_{n}(y) \,x^{n} -a_0(y)\\
 & = & 
\sum_{n=0}^\infty a_{n}(y) \,x^{n+1}
-\sum_{n=0}^\infty b_{n}(y) \,x^{n+1}\\
 & = & x\,T_y(x)-\sum_{n=0}^\infty b_{n}(y) \,x^{n+1}.
\end{eqnarray*}
Using \eqref{4}, we can solve to get
$$
T_y(x)=\frac{{\rm erf}({\sqrt y})-\sum_{n=0}^\infty b_{n}(y) \,x^{n+1}}
{1-x}.
$$
The change of variables $t=x\sin\theta$ gives
$$
{\rm erf}(x)=\frac{2}{\sqrt \pi}\,x\,e^{-x^2}\int_{\theta=0}^{\pi/2}
e^{x^2\cos^2\theta}\cos\theta\,d\theta.
$$
Since the exponential function is entire, we can integrate its
Taylor series and use \eqref{6} to get
\begin{equation}
{\rm erf}(x)  =  \frac{2}{\sqrt \pi}\,e^{-x^2}\sum_{m=0}^\infty
\frac{x^{2m+1}}{m!}\int_{\theta=0}^{\pi/2}
\cos^{2m+1}\!\theta\,d\theta
  =  \frac{2}{\sqrt \pi}\,e^{-x^2}\sum_{m=0}^\infty
\frac{2^m x^{2m+1}}
{1\cdot 3\cdot5\cdots (2m+1)}.\label{7}
\end{equation}
This series converges absolutely for all $x$.

Now, using \eqref{5} and  \eqref{7},
$$
\sum_{n=0}^\infty b_n(y)\,x^{n+1}  =  \frac{e^{-y}}{\sqrt \pi}
\sum_{n=0}^\infty \frac{2^ny^{n+1/2}x^{n+1}}{(n+1/2)1\cdot3\cdots(2n-1)}
  =  {\sqrt x}\,e^{-y} e^{xy}{\rm erf}(\sqrt{xy}).
$$
And, 
$$
T_y(x)=\sum_{n=0}^\infty \frac{x^n\gamma(n+1/2,y)}
{\Gamma(n+1/2)}
=\frac{{\rm erf}({\sqrt y})
-{\sqrt x}\,e^{-y} e^{xy}{\rm erf}(\sqrt{xy})}{1-x}.
$$
This result holds for $0<x<1$ and $y\geq 0$.  If we now put
$x=\beta^2/(1+\beta^2)$ and $y=(1+\beta^2)t$ we obtain the
required identity.

The published solution is completely different, using Laplace transforms \cite{axeness}.

\section{Binomial theorem}\label{sectionbinomial}
Every summation method needs to cut its teeth on the binomial theorem.

\begin{theorem}\label{theorembinomial}
Let $a,z\in\C$ with $\abs{z}<1$.  Then $(1-z)^{-a}=\frac{1}{\Gamma(a)}\sum_{n=0}^\infty
\frac{\Gamma(n+a)z^n}{n!}$.
\end{theorem}
\begin{proof}
First suppose $a\not\in\{0,-1,-2,\ldots\}$.
Write $A_n=\Gamma(n+a)/(\Gamma(a)n!)$ and let $S(z)=\sum_{n=0}^\infty A_nz^n$.  From \eqref{2}, $A_{n+1}/A_n=(n+a)/(n+1)$.
The ratio test then shows the series converges absolutely for $\abs{z}<1$.
It can be differentiated within this circle of convergence, yielding a series with the 
same radius of convergence, both of which are analytic in this circle.

We have the linear recurrence relation $(n+1)A_{n+1}-(n+a)A_n=0$.  Multiply by $z^n$ and
sum from $n=0$ to $\infty$.  This gives
\begin{eqnarray*}
0 & = & \sum_{n=0}^\infty(n+1)A_{n+1}z^n-\sum_{n=0}^\infty(n+a)A_nz^n\\
 & = & \sum_{n=1}^\infty nA_{n}z^{n-1}-z\sum_{n=0}^\infty nA_nz^{n-1}-a\sum_{n=0}^\infty A_nz^n\\
 & = & S'(z) -zS'(z)-aS(z).
\end{eqnarray*}
This gives the differential equation $(1-z)S'(z)-aS(z)=0$ with initial condition
$S(0)=A_0=1$.  The solution is $S(z)=(1-z)^{-a}$.  If $z$ is real the
solution is unique by direct integration of the differential equation.
If $z$ is complex this is still the case since the function $z\mapsto 1/(1-z)$ is
analytic except for the pole at $1$ hence the integration is independent of path.

The function $z\mapsto 1/\Gamma(z)$ is entire, with zeroes at $0$ and negative integers.
If $a=-m\in\{0,-1,-2,\ldots\}$ then, using the reflection identity \cite[6.1.17]{abramowitz}, 
$\Gamma(z)=\pi/(\Gamma(1-z)\sin(\pi z))$,
we get $A_n=\Gamma(m+1)/(\Gamma(m-n+1)n!))=0$ for $n\geq m+1$.  The series then reduces to the
familiar formula $(1-z)^m=\sum_{n=0}^m\binom{m}{n}z^n$, which now holds for all $z\in\C$.
\end{proof}

There are lots of proofs of the binomial theorem in the literature, although many of them
confine $z$ to $\R$ and $-a$ to $\N$.  For example, see \cite{ebiscu}, where a recurrence
relation in $a$ is used.  See \cite{abbas} for a different differential equation method
than what we used and references for various other methods.

\section{Bessel functions}\label{sectionbessel}
The generating function for Bessel functions can be obtained using these methods
with a first order differential equation.

\begin{theorem}\label{theorembesselgf}
For each $z,t\in\C$,
\begin{equation}
\sum_{n=-\infty}^\infty J_n(z)t^n=e^{\frac{z}{2}(t-1/t)}.\label{besselgf}
\end{equation}
\end{theorem}
There are many ways to obtain \eqref{besselgf}, such as Fourier transform and Fourier series 
\cite[\S4.9]{andrewsaskeyroy}.
There are many ways to verify this formula. For example, both $e^{zt/2}$ and $e^{-z/(2t)}$
can be expanded in Taylor series, their product computed, and then powers of $t$ equated.

A few Bessel function identities are required, all of which can be found in \cite{abramowitz}.
\begin{align}
&J_n(z)  \sim  \frac{1}{\sqrt{2\pi n}}\left(\frac{ez}{2n}\right)^n \text{ as } n\to\infty
\label{bessel1}\\
&J_{-n}(z)=(-1)^nJ_n(z)\label{bessel2}\\
&zJ_{n+1}(z)-2nJ_n(z)+zJ_{n-1}(z)=0\label{bessel3}\\
&J_n(0)=\left\{\begin{array}{cl}
1, & n=0\\
0, & n\geq 1
\end{array}
\right.\label{bessel4}\\
&J_n(z)  = J_{n+2}(z)+2J'_{n+1}(z).\label{bessel5}
\end{align}
\begin{proof}
Let $S_z(t)=\sum_{n=-\infty}^\infty J_n(z)t^n$.  By \eqref{bessel1}
and \eqref{bessel2},
the series converges absolutely for each $z,t\in\C$. Similarly with the termwise 
differentiated series 
$S'_z(t)$.  Multiply the recurrence relation \eqref{bessel3} by $t^n$, sum, and re-index:
\begin{eqnarray*}
0 & = & \frac{z}{t}\sum_{n=-\infty}^\infty J_n(z)t^n-2t\sum_{n=-\infty}^\infty nJ_n(z)t^{n-1}
+zt\sum_{n=-\infty}^\infty J_n(z)t^n\\
 & = & \frac{z}{t}S_z(t)-2tS_z'(t)+zt S_z(t).
\end{eqnarray*}
This gives the ordinary differential equation $S_z'(t)-\frac{z}{2}(1+1/t^2)S_z(t)=0$.
The unique solution is $S_z(t)=C(z)e^{(z/2)(t-1/t)}$ for some function $C\fn\C\to\C$.
Use \eqref{bessel5} to get
$$
C(z)=S_z(1)=\sum_{n=-\infty}^\infty J_n(z)=
\sum_{n=-\infty}^\infty (J_n(z)+2J'_n(z))=C(z)+2C'(z).
$$
Therefore, $C'=0$ and $C(z)=C(0)=J_0(0)=1$, from \eqref{bessel2} and \eqref{bessel5}.
\end{proof}

\section{American Mathematical Monthly Problem 11274}\label{section11274}
{\it 
Prove that for nonnegative integers $m$ and $n$,
\begin{equation}
\sum_{k=0}^m 2^k \binom{2m-k}{m+n} =4^m -\sum_{j=1}^n \binom{2m+1}{m+j}.\label{prob}\\
\end{equation}
}
See \cite{11274}.\\

\noindent {\it Solution.}
For nonnegative integers $a$ and $b$ we will use the
identities
\begin{equation}
\binom{a}{b}=\binom{a}{a-b}\label{21}
\end{equation}
and
\begin{equation}
\binom{a}{b} =0\quad\text{for } b>a.\label{22}
\end{equation}

First prove the result for $m<n$.  By \eqref{22}, 
$\sum_{k=0}^m 2^k\binom{2m-k}{m+n}=0$.  And,
$$
\sum_{j=1}^n \binom{2m+1}{m+j}  =  \sum_{j=1}^{m+1} \binom{2m+1}{m+j}
=\sum_{j=0}^m\binom{2m+1}{m+j+1}
  =  \sum_{j=0}^m \binom{2m+1}{j},
$$
using \eqref{21} and then reversing the order of summation.  Now let
$y\in\R$ and write
\begin{eqnarray*}
\sum_{j=0}^m \binom{2m+1}{j}y^j & = & \sum_{j=0}^{2m+1} \binom{2m+1}{j}y^j
- \sum_{j=m+1}^{2m+1} \binom{2m+1}{j}y^j\\
 & = & (y+1)^{2m+1} -\sum_{j=0}^{m} \binom{2m+1}{m+1+j}y^{m+1+j},
\end{eqnarray*}
using the binomial theorem and the change of variables $j\mapsto j+m+1$.
Using \eqref{21} and then reversing the order of summation, we have
\begin{eqnarray}
\sum_{j=0}^m \binom{2m+1}{j}y^j & = & (y+1)^{2m+1} -
\sum_{j=0}^{m} \binom{2m+1}{m-j}y^{m+1+j}\notag\\
 & = & (y+1)^{2m+1} -
\sum_{j=0}^{m} \binom{2m+1}{j}y^{2m+1-j}.\label{result1}
\end{eqnarray}
Put $y=1$ to get 
\begin{equation}
\sum_{j=1}^{m+1}\binom{2m+1}{m+j} = \sum_{j=0}^m\binom{2m+1}{j}=2^{2m}.
\label{26}
\end{equation}
Hence, both sides of \eqref{prob} are zero when $m<n$.

Now suppose $m\geq n$.  We will show that \eqref{prob} is a special
case of  the identity
\begin{equation}
\sum_{k=0}^{m-n}\binom{2m-k}{m+n}x^k = \sum_{j=1}^{m+1}\binom{2m+1}{m+j}(x-1)^{j-n-1}
-\sum_{j=1}^n\binom{2m+1}{m+j}(x-1)^{j-n-1},\label{result2}
\end{equation}
valid for all $x\not=1$.

Consider the left side of \eqref{result2}.
Let $p=m-n$, $q=m+n$ and $k=p-\ell$.  Then
$$
\sum_{k=0}^{m-n}\binom{2m-k}{m+n}x^k = 
x^p\sum_{\ell=0}^p\binom{q+\ell}{q}x^{-\ell}.
$$
Let $A(\ell,q)=\binom{q+\ell}{q}$ and consider $S(p,q;x):=\sum_{\ell=0}^p
A(\ell,q)x^\ell$.  Notice the recurrence relation
$$
(\ell +1)A(\ell+1,q) - (q+\ell+1)A(\ell,q)=0.
$$
Hence, for all $x\in\R$, 
$$
\sum_{\ell=0}^p(\ell +1)A(\ell+1,q)x^\ell -
\sum_{\ell=0}^p (q+\ell+1)A(\ell,q)x^\ell=0.
$$
Suppose $x<1$.  Then, using \eqref{constant} and \eqref{polynomial},
$$
\sum_{\ell=1}^{p+1}\ell\,A(\ell,q)x^{\ell -1} -
(q+1)\sum_{\ell=0}^p A(\ell,q)x^\ell
-x\sum_{\ell=0}^p\ell A(\ell,q)x^{\ell-1}=0.
$$
This then gives the differential equation
$$
S'(p,q;x)+\left(\frac{q+1}{x-1}\right)S(p,q;x)=\left(\frac{p+1}{x-1}\right)
\binom{p+q+1}{q}x^p
$$
with initial condition $S(p,q;0)=\binom{q}{q}=1$.  The unique solution on $(-\infty,1)$ is
$$
S(p,q;x)=(1-x)^{-(q+1)}\left\{1-(p+1)\binom{p+q+1}{q}\int_0^x(1-t)^qt^p\,dt
\right\}.\label{24}
$$
Repeated integration by parts shows that 
\begin{align*}
&\int_0^x(1-t)^qt^p\,dt  =  \int_0^x(1-t)^{m+n}t^{m-n}\,dt\notag\\
&=(m-n)!(m+n)!\left[\frac{1}{(2m+1)!}-\sum_{k=0}^{m-n}
\frac{(1-x)^{m+n+k+1}x^{m-n-k}}{(m+n+k+1)!(m-n-k)!}\right].\label{25}
\end{align*}

This then gives
\begin{eqnarray}
S(p,q;x) & = & S(m-n,m+n;x)\notag\\
 & = & (1-x)^{-(m+n+1)}(2m+1)!
\sum_{k=0}^{m-n}
\frac{(1-x)^{m+n+k+1}x^{m-n-k}}{(m+n+k+1)!(m-n-k)!}\notag\\
 & = & \sum_{k=0}^{m-n}\binom{2m+1}{m+n+k+1}(1-x)^{k}x^{m-n-k}.\label{S}
\end{eqnarray}

Using \eqref{S} with $x>1$, the left side of \eqref{result2}
then becomes
\begin{eqnarray}
\sum_{k=0}^{m-n}\binom{2m-k}{m+n}x^k
 & = &   x^{m-n}S(m-n,m+n;1/x)\label{hypergeometric1}\\
 & = &  \sum_{k=0}^{m-n}\binom{2m+1}{m+n+k+1}(x-1)^k\label{hypergeometric2}\\
 & = &   \sum_{j=n+1}^{m+1}\binom{2m+1}{m+j}(x-1)^{j-n-1}\quad\text{(Let }
j=n+k+1\text{.)}\notag\\
 & = &    \sum_{j=1}^{m+1}\binom{2m+1}{m+j}(x-1)^{j-n-1}
- \sum_{j=1}^{n}\binom{2m+1}{m+j}(x-1)^{j-n-1}.\label{hypergeometric3}
\end{eqnarray}
We then get
\begin{equation}
(x-1)^{n+1}\sum_{k=0}^{m-n}\binom{2m-k}{m+n}x^k =
\sum_{j=1}^{m+1}\binom{2m+1}{m+j}(x-1)^{j} 
-\sum_{j=1}^n\binom{2m+1}{m+j}(x-1)^{j}.\label{Final}
\end{equation}

Putting $x=2$ and using \eqref{26} completes the proof.

Note that \eqref{Final} was proved for $x>1$ but it is a polynomial equality so it
holds for all $x\in\R$.  It therefore also holds everywhere in the complex plane.

This problem was proposed by Donald Knuth
\cite{11274}.  Since my solution proved the
more general results \eqref{result1} and \eqref{Final},
I didn't think it would be chosen for publication
in {\it American Mathematical Monthly}.  (It wasn't.  The published solution was by
Donald Knuth and Julian Hook and was combinatoric in nature \cite{knuthhook}.)
I mailed my solution to Knuth
and he replied with a nice postcard.  

Knuth's web page says, 
{\it Email is a wonderful thing for people whose role in life is to be on top 
of things. But not for me; my role is to be on the bottom of things.}
This explains why he does not use email.

The equality \eqref{hypergeometric1} - \eqref{hypergeometric2}
is a special case of a hypergeometric identity \cite[15.3.6]{abramowitz}.  The hypergeometric function
is ${}_2F_1(a,b;c;x)=\frac{\Gamma(c)}{\Gamma(a)\Gamma(b)}
\sum_{n=0}^\infty\frac{\Gamma(a+n)\Gamma(b+n)x^n}{\Gamma(c+n)n!}$.
We have
\begin{align*}
&\sum_{k=0}^{m-n}\binom{2m-k}{m+n}x^k=\binom{2m}{m-n}{}_2F_1(-(m-n),1;-2m;x)
=\sum_{k=0}^{m-n}\binom{2m+1}{m+n+k+1}(x-1)^k\\
&=\binom{2m+1}{m+n+1}{}_2F_1(-(m-n),1;m+n+2;1-x)
=P_{m-n}^{(m+n+1,-2m-1)}(2x-1).
\end{align*}
The parameter values `$1$' indicate the hypergeometric functions are actually
confluent hypergeometric functions \cite[Chapter~13.]{abramowitz}.
In this case the hypergeometric functions reduce to Jacobi polynomials
\cite[22.5.42]{abramowitz}, for which
the Rodrigues definition is \cite[22.11.1]{abramowitz}
$$
P_n^{(\alpha,\beta)}(x)=\frac{(-1)^n}{2^nn!(1-x)^\alpha(1+x)^\beta}\frac{d^n}{dx^n}
\left[(1-x)^{\alpha+n}(1+x)^{\beta+n}\right].
$$

Furthermore, the first sum in \eqref{hypergeometric3}
can be identified as a Jacobi polynomial.
Using identities as above,
$$
\sum_{j=1}^{m+1}\binom{2m+1}{m+j}(x-1)^{j-n-1} =(x-1)^{-n}P_m^{(m+1,-2m-1)}(2x-1).
$$
Now putting $x=2$ in \eqref{hypergeometric1} - \eqref{hypergeometric3} gives
$$
P_{m-n}^{(m+n+1,-2m-1)}(3)=P_m^{(m+1,-2m-1)}(3)-\sum_{j=1}^n\binom{2m+1}{m+j}.
$$
In this order, the three terms here correspond to the three in the original
statement of the problem \eqref{prob}.  There are several ways of rewriting these
in terms of Jacobi polynomials and confluent hypergeometric functions.  The first
two sums reduced to well known functions since the cut off of the binomial terms
via \eqref{22} agreed with
the upper summation limit.  The third sum is more complicated since
it is an indefinite summation with the arbitrary upper limit $n$.

\end{document}